\documentclass[12pt]{amsart}
 \usepackage{epsfig}
\usepackage{graphicx}
\usepackage{amscd}
\usepackage{amsmath}
\usepackage{amsxtra}
\usepackage{amsfonts}
\usepackage{amssymb}

\newtheorem{theorem}{Theorem}[section]
\newtheorem{corollary}[theorem]{Corollary}

\newtheorem{proposition}[theorem]{Proposition}

\newtheorem{conjecture}[theorem]{Conjecture}
\theoremstyle{definition}
\newtheorem{definition}[theorem]{Definition}
\newtheorem{remark}[theorem]{Remark}

\newtheorem{example}[theorem]{Example}
\theoremstyle{remark}

\renewcommand{\theclaim}{\textup{\theclaim}}

\newtheorem*{acknowledgements}{Acknowledgements}

\numberwithin{equation}{section}

\def\openone

{\mathchoice

{\hbox{\upshape \small1\kern-3.3pt\normalsize1}}

{\hbox{\upshape \small1\kern-3.3pt\normalsize1}}

{\hbox{\upshape \tiny1\kern-2.3pt\SMALL1}}

{\hbox{\upshape \Tiny1\kern-2pt\tiny1}}}

\makeatletter

\newbox\ipbox

\newcommand{\diracb}[1]{\left\langle #1\mathrel{\mathchoice

{\setbox\ipbox=\hbox{$\displaystyle \left\langle\mathstrut
#1\right.$}

\vrule height\ht\ipbox width0.25pt depth\dp\ipbox}

{\setbox\ipbox=\hbox{$\textstyle \left\langle\mathstrut
#1\right.$}

\vrule height\ht\ipbox width0.25pt depth\dp\ipbox}

{\setbox\ipbox=\hbox{$\scriptstyle \left\langle\mathstrut
#1\right.$}

\vrule height\ht\ipbox width0.25pt depth\dp\ipbox}

{\setbox\ipbox=\hbox{$\scriptscriptstyle \left\langle\mathstrut
#1\right.$}

\vrule height\ht\ipbox width0.25pt depth\dp\ipbox}

}\right. }

\newcommand{\dirack}[1]{\left. \mathrel{\mathchoice

{\setbox\ipbox=\hbox{$\displaystyle \left.\mathstrut
#1\right\rangle$}

\vrule height\ht\ipbox width0.25pt depth\dp\ipbox}

{\setbox\ipbox=\hbox{$\textstyle \left.\mathstrut
#1\right\rangle$}

\vrule height\ht\ipbox width0.25pt depth\dp\ipbox}

{\setbox\ipbox=\hbox{$\scriptstyle \left.\mathstrut
#1\right\rangle$}

\vrule height\ht\ipbox width0.25pt depth\dp\ipbox}

{\setbox\ipbox=\hbox{$\scriptscriptstyle \left.\mathstrut
#1\right\rangle$}

\vrule height\ht\ipbox width0.25pt depth\dp\ipbox}

} #1\right\rangle}

\newcommand{\cj}[1]{\overline{#1}}

\newcommand{\bz}{\mathbb{Z}}

\newcommand{\br}{\mathbb{R}}

\newcommand{\bn}{\mathbb{N}}

\newcommand{\vectr}[2]{\left[\begin{array}{r}#1\\#2\end{array}\right]}

\def\blfootnote{\xdef\@thefnmark{}\@footnotetext}

\newcommand{\wdots}{\ldots}
\newcommand{\proj}{\operatorname{proj}}
\renewcommand{\mod}{\operatorname{mod}}

\input xy
\xyoption{all}
\usepackage{amssymb}

\def\-{^{-1}}

\begin{document}
\title[Fourier duality for affine IFSs]{Probability and Fourier duality for affine iterated function systems}
\author{Dorin Ervin Dutkay}
\blfootnote{Research supported in part by a grant from the National Science Foundation DMS-0704191}
\address{[Dorin Ervin Dutkay] University of Central Florida\\
	Department of Mathematics\\
	4000 Central Florida Blvd.\\
	P.O. Box 161364\\
	Orlando, FL 32816-1364\\
U.S.A.\\} \email{ddutkay@mail.ucf.edu}

\author{Palle E.T. Jorgensen}
\address{[Palle E.T. Jorgensen]University of Iowa\\
Department of Mathematics\\
14 MacLean Hall\\
Iowa City, IA 52242-1419\\}\email{jorgen@math.uiowa.edu}
\thanks{} 

\begin{abstract}
   Let $d$ be a positive integer, and let $\mu$ be a finite measure on $\br^d$. In this paper we ask when it is possible to find a subset $\Lambda$ in $\br^d$ such that the corresponding complex exponential functions $e_\lambda$ indexed by $\Lambda$ are orthogonal and total in $L^2(\mu)$. If this happens, we say that $(\mu, \Lambda)$ is a spectral pair. This is a Fourier duality, and the $x$-variable for the $L^2(\mu)$-functions is one side in the duality, while the points in $\Lambda$ is the other. Stated this way, the framework is too wide, and we shall restrict attention to measures $\mu$ which come with an intrinsic scaling symmetry built in and specified by a finite and prescribed system of contractive affine mappings in $\br^d$; an affine iterated function system (IFS). This setting allows us to generate candidates for spectral pairs in such a way that the sets on both sides of the Fourier duality are generated by suitably chosen affine IFSs. For a given affine setup, we spell out the appropriate duality conditions that the two dual IFS-systems must have. Our condition is stated in terms of certain complex Hadamard matrices. Our main results give two ways of building higher dimensional spectral pairs from combinatorial algebra and spectral theory applied to lower dimensional systems.
\end{abstract}
\maketitle \tableofcontents
\section{Introduction}\label{intr}

   The use of traditional Fourier series has up to recently been restricted to the setting of Fourier duality between groups; in the abelian case \cite{Rud62}, between compact groups (such as tori) on the one side, each group coming with its Haar measure; and discrete abelian groups (such as lattices) on the other. However in dynamics and in other applications to computational mathematics, one is often faced with sets arising as attractors, highly non-linear, and coming equipped with equilibrium measures. This has led to attempts at adapting traditional Fourier tools to these non-linear and non-group settings. In this paper we address the Fourier duality question for affine iterated function systems. For some of the earlier literature we refer the reader to \cite{JoPe92,  JoPe93,  JoPe94,  JoPe95,  JoPe98a,  JoPe98b,DuJo07a,DuJo07b, DJ07c, DJ07d, DR07,DHPS08,Jor06,JKS07,OS05}.

 {\it Iterated function systems} (IFS) in $\br^d$ are natural generalizations of more familiar Cantor sets on the real line. Like their linear counterparts, they arise as limit sets $X$ for recursively defined dynamical systems. While the functions used may be affine, the limit $X$ itself will typically be a highly non-linear object, and will include complicated geometries. They arise in operator algebras and in representation theory; and they form models for ``attractors'' in dynamical systems arising in nature. For $d = 2$, the Sierpinski gasket is a notable example, and there is a variety of possibilities for $d>2$ as well. Each affine IFS $X$ possesses (normalized) invariant measures $\mu$ (see \eqref{eq1_4}), naturally associated with the system at hand (denote by $X$ the support of the measure).

     {\bf Question}: When is there some {\it orthogonal basis} in the Hilbert space  $L^2(X,\mu)$ of the form $\{e_\lambda\,|\,\lambda\in\Lambda\}$, where $\Lambda\subset\br^d$ and $e_\lambda(x) = \exp(2\pi i x\cdot\lambda)$?

  {\bf Definition.} The restricted class of IFSs $\mu$ for which a basis $\{e_\lambda\}_{\lambda\in\Lambda}$ can be found are called {\it spectral measures}, and the functions $e_\lambda$ are said to form a {\it Fourier basis}. $(\mu, \Lambda)$ is then called a {\it spectral pair}.

    If $X$ is the middle-third Cantor set, Jorgensen and Pedersen \cite{JoPe98} proved that there is no such Fourier basis. Nonetheless, many spectral measures have since been found within various important classes of fractals; and their significance has been explored by many authors.

   Affine IFSs often take the following form: Start with a $d$ by $d$ matrix $R$ and a finite subset $B$ in $\br^d$. Then consider the associated set of affine mappings in $\br^d$, $\tau_b$ of the form $\tau_b(x) = R^{-1}(x+b)$, where $R$ is further assumed to be an expanding integer matrix and $b\in B\subset\br^d$.

    In \cite{DuJo07b} we conjectured that the {\it spectral measures} arise precisely when there exists {\it duality pairing}, i.e., another system $L$ such that the two define a {\it complex Hadamard matrix} $H$ (see equation \eqref{eqhada}), the order of the matrix $H$ being the cardinality of $B$. We proved the conjecture in \cite{DuJo07b} when an additional assumption, called ``reducibility'', is placed on the triple $(R,B,L)$.

      Two approaches to IFSs have been popular: one based on a discrete version of the more familiar and classical second order Laplace differential operator of {\it potential theory}, see \cite{KSW01, Kig04, LNRG96};  and the second approach is based on {\it Fourier series}, see e.g.,  \cite{JoPe98, DuJo06a}. The first model is motivated by infinite discrete network of resistors, and the harmonic functions are defined by minimizing a global measure of resistance, but this approach does not rely on Fourier series. In contrast, the second approach begins with Fourier series, and it has its classical origins in lacunary Fourier series \cite{Kah86}.

   {\bf Hadamard matrices.}   
Let $R$ be a $d\times d$ integer matrix, $B\subset\bz^d$ and $L\subset\bz^d$ having the same cardinality as $B$, $\#B=\#L=:N$. We call $(R,B,L)$ a {\it Hadamard triple} if 
the matrix
\begin{equation}\label{eqhada}
\frac{1}{\sqrt{N}}(e^{2\pi iR^{-1}b\cdot l})_{b\in B,l\in L}
\end{equation}
is unitary.

Let $\mu_B$ be the invariant measure associated to the affine IFS $\tau_b(x)=R^{-1}(x+b)$, $b\in B$. (See Theorem \ref{th1_1} below).
We conjectured in \cite{DuJo07b} that the existence of a set $L$ such that $(R,B,L)$ is a Hadamard triple is sufficient to obtain orthonormal bases of exponentials in $L^2(\mu_B)$. 
\begin{conjecture}\label{conj}
Let $R$ be a $d\times d$ expansive integer matrix, $B$ a subset of $\bz^d$ with $0\in B$. Let $\mu_B$ be the invariant measure of the associated IFS $(\tau_b)_{b\in B}$. If there exists a subset $L$ of $\bz^d$ such that $(R,B,L)$ is a Hadamard triple and $0\in L$ then $\mu_B$ is a spectral measure. 
\end{conjecture}

In \cite{DuJo07b} we proved that this conjecture is true under a certain ``reducibility'' assumption.

\begin{theorem}
Let $R$ be an expanding $d\times d$ integer matrix, $B$ a subset of $\bz^d$ with $0\in B$. Assume that there exists a subset $L$ of $\bz^d$ with $0\in L$ such that $(R,B,L)$ is a Hadamard triple which satisfies the reducibility condition (see \cite{DuJo07b}). Then the invariant measure $\mu_B$ is a spectral measure. In particular the conjecture is true in dimension $d=1$.
\end{theorem}

In this paper we will show that the conjecture is true also in some other cases, by reducing the problem to some ``building blocks'' in lower dimensions (Theorem \ref{thi1}, Corollary \ref{co2}). We give an example of an affine IFS when the reducibility condition from \cite{DuJo07b} is not satisfied but the associated measure is still spectral (Example \ref{exnr}).
 
\section{Definitions}

The purpose of the present section is to collect the necessary definitions, and to point out the link between the two sides in the general (non-group) Fourier duality; see especially Proposition \ref{prs2} and Definition \ref{def1_2}, Hadamard-triple.
\begin{definition}\label{defs1}
For $\lambda\in\br^d$, let $e_\lambda(x):=e^{2\pi i\lambda\cdot x}$, $x\in\br^d$.

A probability measure on $\br^d$ is called {\it spectral} if there exists a subset $\Lambda$ of $\br^d$ such that the family of exponential functions 
$\{e_\lambda\,|\,\lambda\in\Lambda\}$ is an orthonormal basis for $L^2(\mu)$. In this case $\Lambda$ is called a spectrum for $\mu$.

The Fourier transform of a measure $\mu$ is defined by
\begin{equation}
	\hat\mu(x)=\int e^{2\pi i t\cdot x}\,d\mu(t),\quad(x\in\br^d).
	\label{eqs1}
\end{equation}
\end{definition}

   In a number of applications, one encounters a measure $\mu$ and a subset $\Lambda$ such that the functions $e_\lambda$ indexed by $\Lambda$ are orthogonal in $L^2(\mu)$, but a separate argument is needed in order to show that the family is total. The following is a universal test which applies to any subset $\Lambda$: It is a necessary and sufficient condition on a pair $(\mu, \Lambda)$  allowing us to decide whether $\mu$ is spectral with spectrum $\Lambda$.

\begin{proposition}\label{prs2}\cite{JoPe98,DuJo07b}
Let $\mu$ be a probability measure on $\br^d$. A subset $\Lambda$ of $\br^d$ is a spectrum for $\mu$ iff
\begin{equation}
	\sum_{\lambda\in\Lambda}|\hat\mu(x+\lambda)|^2=1,\quad(x\in\br^d).
	\label{eqs2}
\end{equation}
\end{proposition}
\begin{remark}\label{res1}
Note that in  equation \eqref{eqs2}, we do not have to worry about possible repetitions in $\Lambda$. Indeed, if $\lambda$ appears at least twice in \eqref{eqs2}, then take 
$x=-\lambda$. Since $\mu$ is a probability measure, $\hat\mu(0)=1$. So, the sum on the left of \eqref{eqs2} is at least 2. This would contradict \eqref{eqs2}, and therefore a $\lambda$ can appear at most once in the sum.
\end{remark}

{\bf Affine IFSs.}
Let $R$ be a $d\times d$ expanding integer matrix, i.e., all eigenvalues $\lambda$ satisfy $|\lambda|>1$. Let $B$ be a finite subset of $\bz^d$ of cardinality $\#B=:N$, with $0\in B$. We consider the iterated function system 
\begin{equation}
\tau_b(x)=R^{-1}(x+b),\quad(x\in\br^d,b\in B).	
	\label{eq1_1}
\end{equation}

\begin{theorem}\label{th1_1}\cite{Hut81} There exists a unique compact set $X_B\subset\br^d$ such that
\begin{equation}
	X_B=\bigcup_{b\in B}\tau_b(X_B).
	\label{eq1_2}
\end{equation}
	Moreover 
\begin{equation}
	X_B=\left\{\sum_{k=1}^\infty R^{-k}b_k\,|\, b_k\in B\mbox{ for all }k\in\bn\right\}.
	\label{eq1_3}
\end{equation}
There exists a unique probability measure $\mu=\mu_{R,B}$ on $\br^d$ such that
\begin{equation}
	\int f\,d\mu=\frac1N\sum_{b\in B}\int f\circ\tau_b\,d\mu,
	\label{eq1_4}
\end{equation}
for all continuous functions $f$ on $\br^d$. 
\end{theorem}

{\bf Hadamard triples.}
\begin{definition}\label{def1_2}
Let $R$ and $B$ as above. Let $L\subset\bz^d$. We say that $(R,B,L)$ form a {\it Hadamard triple} if $\#L=\#B=N$, $0\in L$ and the matrix
\begin{equation}
	\frac{1}{\sqrt N}\left( e^{2\pi i R^{-1}b\cdot l}\right)_{b\in B, l\in L}
	\label{eq1_5}
\end{equation}
is unitary.
\end{definition}

{\bf Assumption.} {\it Throughout the paper we will assume that $(R,B,L)$ is a Hadamard triple.}

  An easy computation shows that both Definitions \ref{defs1} and \ref{def1_2} are closed under taking tensor product. For example, the complex Hadamard matrices occurring in \eqref{eq1_5} include the matrices defining the Fourier transform on finite abelian groups, as well as tensor products of these matrices. Currently there is no complete classification of all the complex Hadamard matrices covered by formula \eqref{eq1_5}. Similarly if a fixed complex Hadamard matrix $H$ is given, it is of interest to know all the triples $(R, B, L)$ with the property that $H$ is obtained from $(R, B, L)$ via formula \eqref{eq1_5}. Example \ref{exnr} below serves to illustrate these issues.

     The purpose of this paper is to explore ways of building higher dimensional spectral pairs from combinatorial algebra applied to lower dimensional systems. Theorems \ref{thi1} and Corollary \ref{co2} are cases in point. We also address the converse problem of factoring higher dimensional spectral pairs into products of ``smaller'' ones.

We will denote by $S$ the transpose of $R$, $S:=R^T$. We will need the following ``dual'' iterated function system
\begin{equation}
	\sigma_l(x)=S^{-1}(x+l),\quad(x\in\br^d, l\in L).
	\label{eq1_6}
\end{equation}

\begin{definition}\label{def1_3}
Let 
\begin{equation}
	W_B(x):=\left|\frac{1}{N}\sum_{b\in B}e^{2\pi ib\cdot x}\right|^2,\quad(x\in\br^d).
	\label{eq1_7}
\end{equation}
Then
\begin{equation}
	\sum_{l\in L} W_B(\sigma_l(x))=1.
	\label{eq1_8}
\end{equation}
\end{definition}

\section{Invariant sets and path measures}\label{inse}

 The purpose of the present section is to introduce a setting from symbolic dynamics which we will use on a particular given IFS subject to the Hadamard-triple law (Definition \ref{def1_2}). The data of a Hadamard triple, the scaling matrix  $R$, and the two finite sets $B$, $L$ allow us to create a digital filter $W_B$. We then introduce $W_B$-cycles, and invariant sets. We form a compact space $\Omega$ of infinite words in $L$, and a $W_B$-path-space measure on $\Omega$, see Definition \ref{def1_12}. The purpose of this is to allow us to construct candidates for spectra $\Lambda$, and then to check the condition in Proposition \ref{prs2}.

{\bf Invariant sets.} For $x\in\br^d$, we call a {\it trajectory} of $x$ a set of points $$\{\sigma_{\omega_n}\cdots \sigma_{\omega_1}x\,|\,n\geq1\}$$ where 
$\{\omega_n\}_n$ is a sequence of elements in $L$ such that $W_B(\sigma_{\omega_n}\cdots \sigma_{\omega_1}x)\neq0$ for all $n\geq1$. We denote by $\mathcal{O}(x)$ the union of all trajectories of $x$ and the closure $\cj{\mathcal{O}(x)}$ is called the {\it orbit} of $x$. If $W_B(\sigma_lx)\neq0$ for some $l\in L$ we say that the {\it transition} from $x$ to $\sigma_lx$ is possible.
\par
A closed subset $F\subset\br^d$ is called {\it invariant} if it contains the orbit of all of its points. An invariant subset is called {\it minimal} if it does not contain any proper invariant subsets.
Let 
$$\Omega:=L^{\bn}=\left\{l_1l_2\dots\,|\,l_k\in L,\mbox{ for all }k\in\bn\right\}$$
the space of infinite words over the alphabet $L$.

{\bf Path measures.} For all $x\in\br^d$, there exists a unique probability measure $P_x=P_x(R,B,L)$ on $\Omega$ such that 
for all $l_1,\dots, l_n\in L$, $n\in\bn$
\begin{equation}
	P_x(\{\omega_1\omega_2\dots\in\Omega\,|\,\omega_1=l_1,\dots,\omega_n=l_n\})=\prod_{k=1}^NW_B(\sigma_k\dots\sigma_1x).
	\label{eq1_9}
\end{equation}

The next three propositions can be found in \cite{DuJo07b}.

\begin{proposition}\label{pr1_4}

For all $x\in\br^d$
\begin{equation}
	|\hat\mu(x)|^2=\prod_{k=1}^\infty W_B(S^{-k}x).
	\label{eq1_10}
\end{equation}
For all $x\in\br^d$ and all $l_1l_2\dots\in\Omega$
\begin{equation}
	P_x(\{l_1l_2\dots\})=\prod_{k=1}^\infty W_B(\sigma_{l_k}\dots\sigma_{l_1}x).
	\label{eq1_11}
\end{equation}
\end{proposition}

\begin{proposition}\label{propnf}
Let $F$ be a compact invariant subset. Define 
\begin{equation}
N(F):=\{\omega\in\Omega\,|\, \lim_{n\rightarrow\infty}d(\sigma_{\omega_n}\cdots \sigma_{\omega_1}x,F)=0\}.	
	\label{eq1_12}
\end{equation}
\textup{(}The definition of $N(F)$ does not depend on $x$\/\textup{)}. Define
\begin{equation}
h_F(x):=P_x(N(F)).
	\label{eq1_13}
\end{equation}
Then $0\leq h_F(x)\leq 1$, $h_F$ is continuous,
\begin{equation}
	\sum_{l\in L}W_B(\sigma_l(x))h_F(\sigma_l(x))=h_F(x),\quad(x\in\br^d),
	\label{eqruelle}
\end{equation}

and for $P_x$-a.e.\ $\omega\in\Omega$
\begin{equation}
\lim_{n\rightarrow\infty}h_F(\sigma_{\omega_n}\cdots \sigma_{\omega_1}x)=\left\{\begin{array}{cc}
1,&\mbox{ if }\omega\in N(F),\\
0,&\mbox{ if }\omega\not\in N(F).\end{array}\right.
	\label{eq1_14}
\end{equation}
\end{proposition}

\begin{proposition}\label{propsupp}
Let $F_1, F_2,\dots , F_p$ be a family of mutually disjoint closed invariant subsets of $\br^d$ such that there is no closed invariant set $F$ with $F\cap\bigcup_kF_k=\emptyset$. Then 
$$P_x\left(\bigcup_{k=1}^pN(F_k)\right)=1\quad(x\in\br^d).$$
\end{proposition}

  Recall that the triples $(R, B, L)$ in Definition \ref{def1_2} involve two sets $B$ and $L$. The role they play is that they are the beginning of a Fourier duality based on the scaling with the matrix $R$ on one side of the duality, and with the transposed matrix $S = R^T$ on the other side. The discussion below is based on this, and the tools we build come from the iterated function system based on the pair $(S, L)$ via formulas \eqref{eq1_6}.

\begin{theorem}\label{thcora}\cite[Th\'eor\`eme 2.8]{CCR96}
Let $K$ be minimal compact invariant set contained in the set of zeros of an entire function $h$ on $\br^d$.
\begin{enumerate}
\item[a)] There exists $V$, a proper subspace of $\br^d$ invariant for $S$ \textup{(}possibly reduced to $\{0\}$\textup{)}, such that $K$ is contained in a finite union $\mathcal{R}$ of translates of $V$.
\item[b)] This union contains the translates of $V$ by the elements of a cycle\\  $\{x_0, \sigma_{l_1}x_0,\dots , \sigma_{l_{m-1}}\cdots \sigma_{l_1}x_0 \}$ contained in $K$, and for all $x$ in this cycle, the function $h$ is zero on $x+V$.
\item[c)] Suppose the hypothesis ``(H) modulo $V$'' is satisfied, i.e., for all $p\geq0$ the equality $\sigma_{k_1}\cdots \sigma_{k_p}0-\sigma_{k_1'}\cdots \sigma_{k_p'}0\in V$, with $k_i,k_i'\in L$ implies $k_i-k_i'\in V$ for all $i\in\{1,\dots ,p\}$. Then 
$$\mathcal{R}=\{x_0+V,\sigma_{l_1}x_0+V,\dots ,\sigma_{l_{m-1}}\cdots \sigma_{l_1}x_0+V\},$$
and every possible transition from a point in $K\cap\sigma_{l_q}\cdots \sigma_{l_1}x_0+V$ leads to a point in $K\cap\sigma_{l_{q+1}}\cdots \sigma_{l_1}x_0+V$ for all $1\leq q\leq m-1$, where $\sigma_{l_m}\cdots \sigma_{l_1}x_0=x_0$.
\item[d)] Since the function $W_B$ is entire, the union $\mathcal{R}$ is itself invariant.
\end{enumerate}
\end{theorem}
\par
A particular example of a minimal compact invariant set is a $W_B$-cycle. In this case, the subspace $V$ in Theorem \ref{thcora} can be chosen to be $V=\{0\}$:
\begin{definition} 
A {\it cycle} of length $m$ for the IFS $(\sigma_l)_{l\in L}$ is a set of (distinct) points of the form $\mathcal{C}:=\{x_0,\sigma_{l_1}x_0,\dots ,\sigma_{l_{m-1}}\cdots \sigma_{l_1}x_0\}$, such that $\sigma_{l_m}\cdots \sigma_{l_1}x_0=x_0$, with $l_1,\dots ,l_m\in L$. A $W_B$-{\it cycle} is a cycle $\mathcal{C}$ such that $W_B(x)=1$ for all $x\in\mathcal{C}$.
\par
For a finite sequence $l_1,\dots ,l_m\in L$ we will denote by $\underline{l_1\wdots l_m}$ the path in $\Omega$ obtained by an infinite repetition of this sequence
$$\underline{l_1\wdots l_m}:=(l_1\wdots l_ml_1\wdots l_m\wdots )$$
\end{definition}

\begin{definition}\label{def1_12}
A closed invariant set $F$ is called {\it spectral} if there exists a subset $\Lambda(F)$ of $\br^d$ such that 
\begin{equation}
	P_x(N(F))=\sum_{\lambda\in\Lambda(F)}|\hat\mu(x+\lambda)|^2,\quad(x\in\br^d).
	\label{eq1_17}
\end{equation}
In this case $\Lambda(F)$ is called a {\it spectrum} for $F$.
\end{definition}

\begin{proposition}\label{pr1_13}
Suppose $(F_i)_{i=1,n}$ are invariant sets with the following properties: 
\begin{enumerate}
\item The set $F_i$ is spectral with spectrum $\Lambda(F_i)$;
\item
For each $1\leq k\leq n$ and $1\leq i_1<i_2<\dots<i_k\leq n$, the set 
$F_{i_1}\cap F_{i_2}\cap\dots\cap F_{i_k}$ is spectral, with spectrum $\Lambda(F_{i_1})\cap\dots\cap\Lambda(F_{i_k})$.

Then 
\begin{equation}
	P_x(\bigcup_{i=1}^n N(F_i))=\sum_{\lambda\in\cup_{i=1}^n\Lambda(F_i)}|\hat\mu(x+\lambda)|^2
	\label{eq1_18}
\end{equation}
\end{enumerate}
\end{proposition}

\begin{proof}
First we prove that if $F_1$ and $F_2$ are invariant sets then
\begin{equation}
	N(F_1\cap F_2)=N(F_1)\cap N(F_2).
	\label{eq1_19}
\end{equation}
The inclusion ``$\subset$'' is clear. For the converse, let $\omega\in N(F_1)\cap\ N(F_2)$. Then
$$d(\sigma_{\omega_n}\dots\sigma_{\omega_1}(0),F_i)\rightarrow\, 0\mbox{ for }i=1,2.$$ Suppose, by contradiction, that $d(\sigma_{\omega_n}\dots\sigma_{\omega_1},F_1\cap F_2)$ does not converge to 0. Then there exists a $\delta>0$ and a subsequence such that
$y_n:=\sigma_{\omega_{k(n)}}\dots\sigma_{\omega_1}0$ and $d(y_n,F_1\cap F_2)\geq\delta$. Since $0$ is contained in the attractor $X_L$ of the IFS $(\sigma_l)_{l\in L}$, it follows that $y_n\in X_L$ for all $n$. Since $X_L$ is compact we can find a subsequence such that
$\sigma_{\omega_{k(n(p))}}\dots\sigma_{\omega_1}0$ converges to some point $x_0$. But then $x_0$ must be both in $F_1$ and $F_2$, which gives the contradiction.

By induction, we can extend \eqref{eq1_19} to any finite number of invariant sets. (Note that the intersection of invariant sets is itself invariant)

We have then
$$P_x(\bigcup_{i=1}^n N(F_i))=\sum_{k=1}^n(-1)^{k+1}\sum_{1\leq i_1<\dots <i_k\leq n} P_x(\bigcap_{j=1}^kN(F_{i_j}))=$$$$
\sum_{k=1}^n(-1)^{k+1}\sum_{1\leq i_1<\dots <i_k\leq n} P_x(N(\bigcap_{j=1}^kF_{i_j}))=
\sum_{k=1}^n(-1)^{k+1}\sum_{1\leq i_1<\dots <i_k\leq n}\sum_{\lambda\in\Lambda(\cap_{j=1}^k F_{i_j})}|\hat\mu(x+\lambda)|^2=$$$$
\sum_{k=1}^n(-1)^{k+1}\sum_{1\leq i_1<\dots <i_k\leq n}\sum_{\lambda\in\cap_{j=1}^k\Lambda(F_{i_j})}|\hat\mu(x+\lambda)|^2=
\sum_{\lambda\in\cup_{i=1}^n\Lambda(F_i)}|\hat\mu(x+\lambda)|^2.$$
\end{proof}

Starting with the iterated function system from the pair $(S, L)$ in formulas \eqref{eq1_6}, we have built a system of $W_B$-cycles $\mathcal C$. The proposition below gives conditions for being able to generate a spectrum $\Lambda(\mathcal C)$ directly from this data.

\begin{proposition}\label{pr1_14}
Let $\mathcal C$ be a $W_B$-cycle. Then $\mathcal C$ is a spectral invariant set with spectrum $\Lambda(\mathcal C)=$ the smallest set $\Lambda$ that contains $-\mathcal C$ with the property that $S\Lambda+L\subset\Lambda$.
\end{proposition}

\begin{proof}
Let $\mathcal C=:\{x_0,x_1,\dots,x_{m-1}\}$ with $x_k=\sigma_{l_k}x_{k-1}$ for all $k\in\{2,\dots,m\}$ and $\sigma_{l_m}x_{m-1}=x_0$. 
With \cite[Lemma 4.1]{DuJo07b}, we know that 
$N(\mathcal C)$ consists of the infinite words in $\Omega$ that end in a repetition of the finite word $l_1\dots l_m$, i.e., words of the form 
$\omega_0\omega_1\dots\omega_n\underline{l_1\dots l_m}$.

We proved in \cite[Lemma 4.9]{DuJo07b} that if $\omega=\omega_0\dots\omega_{km-1}\underline{l_1\dots l_m}$, for some $k\geq 0$, $\omega_i\in L$, then
\begin{equation}
	P_x(\{\omega\})=|\hat\mu(x+k_{\mathcal C}(\omega))|^2,\quad(x\in\br^d)|
	\label{eq1_20}
\end{equation}
where
\begin{equation}
	k_{\mathcal C}(\omega)=\omega_0+S\omega_1+\dots+S^{km-1}\omega_{km-1}-S^{km}x_0.
	\label{eq1_21}
\end{equation}

Note that every $\omega\in N(\mathcal C)$ is of the form $\omega_0\dots\omega_{km-1}\underline{\tilde l_1\dots\tilde l_m}$, where 
$\tilde l_1\dots\tilde l_m$ is a circular permutation of $l_1\dots l_m$.

Using \eqref{eq1_20}, we have to prove only that the set $\Lambda(k_{\mathcal C})$ of all numbers
$k_{\mathcal C}(\omega)$ with $\omega$ of the form $\omega_0\dots\omega_{km-1}\underline{\tilde l_1\dots\tilde l_m}$, is the set $\Lambda(\mathcal C)$ defined in the statement of the Proposition.

Since $\sigma_{l_i}x_{i-1}=x_{i}$, we have $-x_{i-1}=l_i-Sx_i$ for all $i$.

First, note that if $\omega=\omega_0\dots\omega_{km-1}\underline{l_1\dots l_m}$, and $\omega_{-1}\in L$, then
$$Sk_{\mathcal C}(\omega)+\omega_{-1}=\omega_{-1}+S\omega_0+\dots+\dots+S^{km}\omega_{km-1}-S^{km+1}x_0=$$$$
\omega_{-1}+\dots+S^{km}\omega_{km}+S^{km+1}l_1-S^{km+2}x_{1}=$$$$\dots=
\omega_{-1}+\dots+S^{km}\omega_{km}+S^{km+1}l_1+\dots+S^{km+m-1}l_{m-1}-S^{km+m}x_{m-1}=$$$$k_{\mathcal C}(\omega_{-1}\dots\omega_{km-1}l_1\dots l_{m-1}\underline{l_ml_1\dots l_{m-1}})$$
This shows that $S\Lambda( k_{\mathcal C})+L\subset\Lambda(k_{\mathcal C})$. On the other hand this calculation shows that any $k_C(\omega)$ from the set $\Lambda(k_{\mathcal C})$ can be obtained (in a unique way!) from a point in the cycle by applying operations of the form $x\mapsto Sx+l$ with $l\in L$. 
This implies that $\Lambda( k_{\mathcal C})=\Lambda(\mathcal C)$.

The uniqueness of the operations $x\mapsto Sx+l$ that lead from a cycle point to a $k_{\mathcal C}(\omega)$ comes from the fact that elements in $L$ are incongruent $\mod S\bz^d$ (see \cite{DuJo06b,DuJo07b} for more details).
\end{proof}

The next result shows how the spectral problem transforms under coordinate changes.

\begin{definition}
We say that two Hadamard triples $(R_1,B_1,L_1)$ and $(R_2,B_2,L_2)$ are {\it conjugate} if there exists a matrix $M\in GL_d(\bz)$ (i.e., $M$ is invertible, and $M$ and $M^{-1}$ have integer entries) such that $R_2=MR_1M^{-1}$, $B_2=MB_1$ and $L_2=(M^T)^{-1}L_1$.
\end{definition}

If the two systems are conjugate then the transition between the IFSs $(\tau_b)_{b\in B_1}$ and $(\tau_{Mb})_{b\in B_1}$ is done by the matrix $M$; and the transition between the IFSs $(\sigma_l)_{l\in L_1}$ and $(\sigma_{(M^T)^{-1}l})_{l\in L_1}$ is done by the matrix $(M^T)^{-1}$.
\begin{proposition}\label{prop2_8}
If $(R_1,B_1,L_1)$ and $(R_2,B_2,L_2)$ are conjugate through the matrix $M$, then
\begin{enumerate}
\item\label{prop2_8(1)}
$\tau_{Mb_1}(Mx)=M\tau_{b_1}(x)$, $\sigma_{(M^T)^{-1}l_1}((M^T)^{-1}x)=(M^T)^{-1}\sigma_{l_1}(x)$, for all $b_1\in B_1$, $l_1\in L_1$;
\item\label{prop2_8(2)}
$W_{B_2}((M^T)^{-1}x)=W_{B_1}(x)$ for all $x\in\br^d$;

\item\label{prop2_8(3)}
For the Fourier transform of the corresponding invariant measures, the following relation holds:
$\hat\mu_{B_2}((M^T)^{-1}x)=\hat\mu_{B_1}(x)$ for all $x\in\br^d$;

\item\label{prop2_8(4)}
For $\omega\in\Omega_1=L_1^{\bn}$, let $(M^T)^{-1}\omega:=((M^T)^{-1}\omega_1,(M^T)^{-1}\omega_2,\dots)$.
The associated path measures satisfy the following relation:
$$P_{(M^T)^{-1}x}^2((M^T)^{-1}E)=P_{x}^1(E).$$
\item\label{prop2_8(5)}
If $F_1$ is invariant for the first IFS, then $(M^T)^{-1}F_1$ is invariant for the second IFS. If in addition $F_1$ has spectrum $\Lambda(F_1)$, then 
$(M^T)^{-1}F_1$ has spectrum $(M^T)^{-1}\Lambda(F_1)$. 
\end{enumerate}
\end{proposition}

\begin{proof}
The results follow by direct computation.
\end{proof}

\section{Invariant subspaces}\label{inva}

 In this section we give two ways of building higher dimensional spectral pairs from combinatorial algebra and spectral theory applied to lower dimensional systems. Theorems \ref{thi1} and Corollary \ref{co2} are cases in point.
 
 \begin{remark}
 If $V$ is a {\it rational} invariant subspace for $S$, i.e., it has a basis consisting of vectors with rational components then, by \cite[Lemme 4.2]{CCR96}, there exists a matrix $M\in GL_d(\bz)$ as in Proposition \ref{prop2_8} that maps $V$ into $\br^r\times\{0\}$, with $r=\dim V$. Therefore we can reduce the study to the case when $V=\br^r\times\{0\}$.
 \end{remark}

Suppose the subspace $\br^r\times\{0\}$ is invariant for $S$.
Then the matrix $S$ has the form 
\begin{equation}
S=\begin{bmatrix}
S_1&C\\
0& S_2\end{bmatrix}
\label{eq2_1}
\end{equation}

where $S_1$ is $r\times r$, $S_2$ is $(d-r)\times(d-r)$, $C$ is $(d-r)\times r$ and $0$ is $r\times (d-r)$.

The matrix $R$ has the form:
\begin{equation}
R=\left[\begin{array}{cc}
A_1&0\\
C^*&A_2\end{array}\right],\mbox{ and }
R^{-1}=\left[\begin{array}{cc}
A_1^{-1}&0\\
-A_2^{-1}C^*A_1^{-1}&A_2^{-1}\end{array}\right],
\label{eq2_2}
\end{equation}

with $A_1=S_1^T$, $A_2=S_2^T$, $C^*=C^T$.
By induction,
\begin{equation}
R^{-k}=\left[\begin{array}{cc}
A_1^{-k}&0\\
D_k&A_2^{-k}\end{array}\right],\mbox{ where }D_k:=-\sum_{l=0}^{k-1}A_2^{-(l+1)}C^*A_1^{-(k-l)}.
\label{eq2_3}
\end{equation}

The set $B$ can be written as 
\begin{equation}
B=\{(r_i,\eta_{i,j})\,|\, i\in\{1,\dots,N_1\}, j\in\{1,\dots, N_2(i)\}\}
\label{eq2_4}
\end{equation}
with $\{r_1,\dots, r_{N_1}\}=\proj_{\br^r}B$.

We have 
\[
X_B=\{\sum_{k=1}^\infty R^{-k}b_k\,|\,b_k\in B\}.
\]
Therefore any element $(x,y)$ in $X_B$ can be written in the following form:
\[
x=\sum_{k=1}^\infty A_1^{-k}r_{i_k},\quad y=\sum_{k=1}^\infty D_kr_{i_k}+\sum_{k=1}^\infty A_2^{-k}\eta_{i_k,j_k}.
\]
\par

Define 
\[
X_1:=\{\sum_{k=1}^\infty A_1^{-k}r_{i_k}\,|\,i_k\in\{1,\dots ,N_1\}\}.
\]
Let $\mu_1$ be the invariant measure for the iterated function system 
\[
\tau_{r_i}(x)=A_1^{-1}(x+r_i),\quad i\in\{1,\dots ,N_1\}.
\]

The set $X_1$ is the attractor of this iterated function system.

{\bf Assumption.} {\it We will assume that the measure $\mu_1$ has no overlap, i.e., $\mu_1(\tau_{r_i}(X_1)\cap\tau_{r_j}(X_1))=0$ for $i\neq j$.}

We will also use the contractions
$$\sigma_l(x)=S_2^{-1}(x+l),\quad(x,l\in\br^{d-r}).$$
\par
For each sequence $\omega=(i_1i_2\wdots )\in\{1,\dots ,N_1\}^{\bn}=:\Omega_1$, define $x(\omega)=\sum_{k=1}^\infty A_1^{-k}r_{i_k}$. Also, because of the non-overlap condition, for $\mu_1$-a.e.\ $x\in X_1$, there is a {\it unique} $\omega$ such that 
$x(\omega)=x$. We define this as $\omega(x)$. This establishes an a.e.\ bijective correspondence between $\Omega_1$ and $X_1$, $\omega\leftrightarrow x(\omega)$.
\par
For $\omega=(i_1i_2\wdots )\in\Omega_1$ define
\[
\Omega_2(\omega):=\{\eta_{i_1,j_1}\eta_{i_2,j_2}\wdots \eta_{i_n,j_n}\wdots \,|\,j_k\in\{1,\dots ,N_2(i_k)\}\}.
\]
For $\omega\in\Omega_1$
\begin{equation}
g(\omega):=\sum_{k=1}^\infty D_kr_{i_k},\mbox{ and }g(x):=g(\omega(x)).
\label{eq2_7}
\end{equation}
Also we denote 
$\Omega_2(x):=\Omega_2(\omega(x))$.
\par
For $x\in X_1$, define
\[
X_2(x):=X_2(\omega(x)):=
\left\{\sum_{k=1}^\infty A_2^{-k}\eta_{i_k,j_k}\biggm|j_k\in\{1,\dots ,N_2(i_k)\}\mbox{ for all }k\right\}.
\]
\par

Note that the attractor $X_B$ has the following form:
\[
X_B=\{(x,g(x)+y)\,|\,x\in X_1,y\in X_2(x)\}.
\]
We will show that the measure $\mu$ can also be decomposed as a product between the measure $\mu_1$ and some measures $\mu_{\omega}^2$ on $X_2(\omega)$.
\par
On $\Omega_2(\omega)$, consider the product probability measure $\mu(\omega)$ which assigns to each $\eta_{i_k,j_k}$ equal probabilities $1/N_2(i_k)$. 
\par
Next we define the measure $\mu_{\omega}^2$ on $X_2(\omega)$. Let 
$r_\omega:\Omega_2(\omega)\rightarrow X_2(\omega)$,
\[
r_\omega(\eta_{i_1,j_1}\eta_{i_2,j_2}\wdots )=\sum_{k=1}^\infty A_2^{-k}\eta_{i_k,j_k}.
\]
\par

Define the measure $\mu_x^2:=\mu_{\omega(x)}^2:=\mu{\omega(x)}\circ r_{\omega(x)}^{-1}$.

  In building higher dimensional spectral pairs from lower dimensional systems by factoring it turns out that overlap properties in the lower dimensional systems play a critical role. This is made precise in the next Proposition. The proof is contained in \cite{DuJo07b}; note that the more restrictive assumptions used there are not needed for this Proposition.

\begin{proposition}\label{pri1}
Using the notations above, suppose 
\begin{enumerate}
\item[(a)] The measure $\mu_1$ has no overlap, i.e., $\mu_1(\tau_{r_i}(X_1)\cap\tau_{r_j}(X_1))=0$ for all $i\neq j$.
\item[(b)] The numbers $N_2(i)=\{\eta\,|\,(r_i,\eta)\in B\}$ are all the same, $N_2(i)=N_2$, $i\in\{1,\dots,N_1\}$. In particular $N_1N_2=N$.
\end{enumerate}
Let $\sigma$ be the shift on $\Omega_1$, $\sigma(i_1i_2\wdots )=(i_2i_3\wdots )$. Let $\omega=(i_1i_2\wdots )\in\Omega_1$. 
\begin{enumerate}
\item
For all measurable sets $E$ in $X_2(\omega)$,
\[
\mu_\omega^2(E)=\frac{1}{N_2}\sum_{j=1}^{N_2(i_1)}\mu_{\sigma(\omega)}^2(\tau_{\eta_{i_1,j}}^{-1}(E)).
\]
The Fourier transform of the measure $\mu_\omega^2$ satisfies the equation:
\begin{equation}\label{eqmu2}
\hat\mu_\omega^2(y)=m(S_2^{-1}y,i_1)\hat\mu_{\sigma(\omega)}^2(S_2^{-1}y),
\end{equation}
where
$$m(y,i_1)=\frac{1}{N_2}\sum_{j=1}^{N_2}e^{2\pi i\eta_{i_1,j}\cdot y}.$$
\item
For all continuous functions on $\br^d$:
\[
\int_{X_B}f\,d\mu=\int_{X_1}\int_{X_2(x)}f(x,y+g(x))\,d\mu_x^2(y)\,d\mu_1(x).
\]
\item
If $\Lambda_1$ is a spectrum for the measure $\mu_1$, then 
$$F(y):=\sum_{\lambda_1\in\Lambda_1}|\hat\mu(x+\lambda_1,y)|^2=\int_{X_1}|\hat\mu_s^2(y)|^2\,d\mu_1(s)\quad(x\in\br^r,y\in\br^{d-r}).$$
Moreover, if 
\begin{equation}
	\tilde W(y):=\frac{1}{N_1}\sum_{i=1}^{N_1}|m(y,i)|^2,\quad(y\in\br^{d-r}),
	\label{eqtw}
\end{equation}
\end{enumerate}
then
\begin{equation}
	F(y)=\tilde W(S_2^{-1}y)F(S_2^{-1}y),\mbox{ and } F(y)=\prod_{k=1}^\infty\tilde W(S_2^{-k}y),\quad(y\in\br^{d-r}).
	\label{eqF}
\end{equation}
\end{proposition}
\begin{definition}
We say that $p\in\br^d$ is an {\it $S$-period} for a function $f$ on $\br^d$ if 
$$f(x+S^np)=f(x),\mbox{ for all }x\in\br^d, n\geq0.$$
\end{definition}

\begin{theorem}\label{thi1}
Suppose $y_0\in\br^{d-r}$ and the set $\br^r\times\{y_0\}$ is invariant. Then $\br^r\times\{0\}$ is invariant for $S$. Using the notations above, there exists $(l_1^0,l_2^0)\in L$ such that $\sigma_{l_2^0}(y_0)=y_0$. 

Assume that the conditions (a) and (b) in Proposition \ref{pri1} are satisfied. Assume in addition that the following conditions are satisfied: 
\begin{enumerate}
\item
The measure $\mu_1$ is spectral with spectrum $\Lambda_1$.
\item The spectrum $\Lambda_1$ has the following properties:
\begin{equation}
	S(\Lambda\times\{-y_0\})+L\supset\Lambda\times\{-y_0\}.
	\label{eqi1_1}
\end{equation}
and
\begin{equation}
	(\lambda_1,-y_0)\mbox{ is an }S\mbox{-period for }W_B,\quad(\lambda_1\in\Lambda_1).
	\label{eqi1_2}
\end{equation}
\item	$(A_1,\{r_i\,|\,i\in\{1,\dots,N_1\},L_1(l_2^0):=\{l_1\,|\,(l_1,l_2^0)\in L\})$ is a Hadamard triple.
\end{enumerate}	
Then $\br^r\times\{y_0\}$ is a spectral set with spectrum
\begin{equation}
\Lambda:=\bigcup_{n\in\bn}(L+SL+\dots+S^{n-1}L+S^n(\Lambda_1\times\{-y_0\})).	
	\label{eqi1_3}
\end{equation}

\end{theorem}

\begin{proof}
Take $x\in\br^r$, then from \eqref{eq1_8} we deduce that there exists $l=(l_1^0,l_2^0)\in L$ such that $W_B(\sigma_l(x,y_0))\neq 0$. Then for $x'\in\br^r$ small enough, we have that $W_B(\sigma_l(x+x',y_0))\neq 0$. Therefore the transition $(x+x',y_0)\mapsto \sigma_l(x+x',y_0)$ is possible. But, since $\br^r\times\{y_0\}$ is invariant, it follows that $\br^r\times\{y_0\}\ni\sigma_l(x+x',y_0)=S^{-1}(x+x',y_0)$. Subtracting, we get
$S^{-1}(x',0)\in\br^{r}\times\{0\}$. Then we can multiply by scalars to see that $S^{-1}(x,0)\in\br^r\times\{0\}$ for all $x\in\br^r$, so $\br^r\times\{0\}$ is invariant for $S$. 

Then we have $$\sigma_l(x,y_0)=S^{-1}(x+l_1^0,y_0+l_2^0)=(S_1^{-1}(x+l_1^0)+ D_1^T(y_0+l_2^0),S_2^{-1}(y_0+l_2^0))\in\br^r\times\{y_0\}$$
which implies that $\sigma_{l_2^0}(y_0)=y_0$.

Next, we compute for $(x,y)\in\br^d$:
\begin{multline*}
\sum_{l_1\in L_1(l_2^0)}W_B(\sigma_{(l_1,l_2^0)}(x,y))=\sum_{l_1}\frac{1}{N^2}\sum_{i,i'=1}^{N_1}\sum_{j,j'=1}^{N_2}e^{2\pi i(r_i-r_i')\cdot(S_1^{-1}(x+l_1)+D_1^T(y+l_2^0))+(\eta_{i,j}-\eta_{i',j'})\cdot(S_2^{-1}(y+l_2^0))}
\\ =(\ast)
\end{multline*}
But using the Hadamard property in (iii) we have
$$\frac{1}{N_1}\sum_{l_1\in L_1(l_2^0)}e^{2\pi i(r_i-r_i')\cdot S_1^{-1}l_1}=\left\{\begin{array}{cc}
1,&i=i'\\
0,&i\neq i'.
\end{array}\right.$$
So 
$$(\ast)=\frac{1}{N_1N_2^2}\sum_{i}\sum_{j,j'}e^{2\pi i(\eta_{i,j}-\eta_{i',j'})\cdot S_2^{-1}(y+l_2^0)}=\frac{1}{N_1}\sum_{i=1}^{N_1}|m(\sigma_{l_2^0}y,i)|^2=\tilde W(\sigma_{l_2^0}(y)).$$
We obtain the equation
\begin{equation}
	\sum_{l_1\in L_1(l_2^0)}W_B(\sigma_{(l_1,l_2^0)}(x,y))=\tilde W(\sigma_{l_2^0}(y)),\quad(y\in\br^{d-r}).
	\label{eq2_2_1}
\end{equation}

Using \cite[Lemma 4.1]{DuJo07b}, we have that $N(\br^r\times\{y_0\})$ consists of all the infinite words in $L^{\bn}$ for which the second components are $l_2^0$ from some point on.

We show first that
\begin{equation}
P_{(x,y)}(\{(\omega_1\wdots \omega_n\wdots )\,|\,\omega_{n,2}=l_{2}^0\mbox{ for all }n\})=\prod_{k=1}^\infty\tilde W(\underbrace{\sigma_{l_{2}^0}\cdots \sigma_{l_{2}^0}}_{k\mbox{ times }}y).	
	\label{eq2_2_2}
\end{equation}

We compute for all $n$, by summing over all the possibilities for the first component, and using (\ref{eq1_11}):
\begin{multline*}
P_{(x,y)}(\{(\omega_1\omega_2\wdots )\,|\,\omega_{k,2}=l_{2}^0, 1\leq k\leq n\})
=
\sum_{l_{1,1},\dots ,l_{n,1}}\prod_{k=1}^nW_B(\sigma_{(l_{k,1},l_{2}^0)}\cdots \sigma_{(l_{1,1},l_{2}^0)}(x,y))
=(*).
\end{multline*}
Using \eqref{eq2_2_1} we obtain further
$$
(*)=\tilde W(\sigma_{l_{2}^0}\cdots \sigma_{l_{2}^0}y)\sum_{l_{1,1},\dots ,l_{n-1,1}}\prod_{k=1}^{n-1}W_B(\sigma_{(l_{k,1},l_{2}^0)}\cdots \sigma_{(l_{1,1},l_{2}^0)}(x,y))
=\dots $$$$
=\prod_{k=1}^n\tilde W(\underbrace{\sigma_{l_{2}^0}\cdots \sigma_{l_{2}^0}}_{k\mbox{ times }}y).
$$
Then, letting $n\rightarrow\infty$ we obtain \eqref{eq2_2_2}.

We have 
$$\sigma_{l_2^0}(y)=S_2^{-1}(y+l_2^0)=S_2^{-1}(y-y_0)+S_2^{-1}(y_0+l_2^0)=S_2^{-1}(y-y_0)+y_0.$$
$$\sigma_{l_2^0}\sigma_{l_2^0}(y)=S_2^{-2}(y-y_0)+S_2^{-1}(y_0+l_2^0)=S_2^{-2}(y-y_0)+y_0.$$
By induction
\begin{equation}
	\underbrace{\sigma_{l_{2}^0}\cdots \sigma_{l_{2}^0}}_{k\mbox{ times }}=S_2^{-k}(y-y_0)+y_0,\quad(k\geq1)
	\label{eq2_2_3}
\end{equation}
Therefore
\begin{equation}
	P_{(x,y)}(\{(\omega_1\wdots \omega_n\wdots )\,|\,\omega_{n,2}=l_{2}^0\mbox{ for all }n\})=\prod_{k=1}^\infty\tilde W(S_2^{-k}(y-y_0)+y_0).
	\label{eq2_2_4}
\end{equation}

Since $\br^r\times\{y_0\}$ is invariant, we have $W_B(\sigma_l(x,y_0))=0$ for all $x\in\br^r$ and $l=(l_1,l_2)\in L$ with $l_2\neq l_2^0$. This implies
that 
$$\tilde W(y_0)=\tilde W(\sigma_{l_2^0}y_0)\stackrel{\mbox{by \eqref{eq2_2_1}}}{=}\sum_{l_1\in L_1(l_2^0)}W_B(\sigma_{(l_1,l_2^0)}(0,y_0))=\sum_{l\in L}W_B(\sigma_l(0,y_0))\stackrel{\mbox{by \eqref{eq1_8}}}{=}1.$$
But then, with \eqref{eqtw}, $|m(y_0,i)|=1$ for all $i\in\{1,\dots,N_1\}$. Using the formula for $m(y_0,i)$ we get that for each $i\in\{1,\dots,N_1\}$ the numbers 
$e^{2\pi i\eta_{i,j}}$, $j\in\{1,\dots,N_2\}$ are all the same (otherwise $|m(y_0,i)|<1$). And this shows that 
$|m(y+y_0,i)|=|m(y,i)|$ so $\tilde W(y+y_0)=\tilde W(y)$ for all $y\in\br^{d-r}$. 

Plug this into \eqref{eq2_2_4}, we get using Proposition \ref{pri1}(iii) and the notations there:
\begin{multline}
	P_{(x,y)}(\{(\omega_1\wdots \omega_n\wdots )\,|\,\omega_{n,2}=l_{2}^0\mbox{ for all }n\})=\prod_{k=1}^\infty\tilde W(S_2^{-k}(y-y_0))=F(y-y_0)\\=\sum_{\lambda_1\in\Lambda_1}|\hat\mu(x+\lambda_1,y-y_0)|^2.
	\label{eq2_2_5}
\end{multline}

Next, fix $l_1,\dots,l_n\in L$. We compute, using \eqref{eq1_11} and \eqref{eq2_2_4}
\begin{multline}
	P_{(x,y)}(\{\omega\,|\,\omega_1=l_1,\dots,\omega_n=l_n,\omega_{k,2}=l_2^0\mbox{ for }k\geq n+1\})=\\
	W_B(\sigma_{l_1}(x,y))\dots W_B(\sigma_{l_n}\dots\sigma_{l_1}(x,y))P_{\sigma_{l_n}\dots\sigma_{l_1}(x,y)}(\{\omega\,|\,\omega_{k,2}=l_2^0,\mbox{ for all }k\geq 1\})=\\
	W_B(\sigma_{l_1}(x,y))\dots W_B(\sigma_{l_n}\dots\sigma_{l_1}(x,y))\sum_{\lambda_1\in\Lambda_1}|\hat\mu(\sigma_{l_n}\dots\sigma_{l_1}(x,y)+(\lambda_1,-y_0))|^2.
	\label{eq2_2_6}
\end{multline}

But 
$$W_B(\sigma_l(x,y))|\hat\mu(\sigma_l(x,y)+(\lambda_1,-y_0))|^2= W_B(S^{-1}((x,y)+l))|\hat\mu(S^{-1}((x,y)+l+S(\lambda_1,-y_0)))|^2$$$$
\stackrel{\mbox{ by \eqref{eqi1_2}}}{=}
W_B(S^{-1}((x,y)+l+S(\lambda_1,-y_0)))|\hat\mu(S^{-1}((x,y)+l+S(\lambda_1,-y_0)))|^2$$$$=|\hat\mu((x,y)+l+S(\lambda_1,-y_0))|^2.$$
By induction,
\begin{multline*}
W_B(\sigma_{l_1}(x,y))\dots W_B(\sigma_{l_n}\dots\sigma_{l_1}(x,y))|\hat\mu(\sigma_{l_n}\dots\sigma_{l_1}(x,y)+(\lambda_1,-y_0))|^2=\\
|\hat\mu((x,y)+l_n+Sl_{n-1}+\dots+S^{n-1}l_1+S^n(\lambda_1,-y_0))|^2.
\end{multline*}
Using this in \eqref{eq2_2_6}, we get
\begin{multline}
	P_{(x,y)}(\{\omega\,|\,\omega_1\in L,\dots,\omega_n\in L,\omega_{k,2}=l_{2,0}\mbox{ for }k\geq n+1\})=\\
	\sum_{l_1,\dots,l_n\in L,\lambda_1\in\Lambda}|\hat\mu((x,y)+l_n+Sl_{n-1}+\dots+S^{n-1}l_1+S^n(\lambda_1,-y_0))|^2.
	\label{eq2_2_7}
\end{multline}

We know by \cite[Lemma 4.1]{DuJo07b}, that $N(\br^r\times\{y_0\})$ consists of all the infinite words in $L^{\bn}$ for which the second component are $l_2^0$ from some point on.

By \eqref{eqi1_1}, the sets $L+SL+\dots+S^{n-1}L+S^n(\Lambda_1\times\{-y_0\})$ are increasing with $n$.  Therefore
$$P_{(x,y)}(N(\br^r\times\{y_0\}))=\sum_{\lambda\in\cup_n(L+SL+\dots+S^{n-1}L+S^n(\Lambda_1\times\{-y_0\}))}|\hat\mu(x+\lambda)|^2.$$

\end{proof}

In the next result we combine the conditions from the theorem into a form that can be used in checking concrete examples.

\begin{corollary}\label{co2}
Suppose the set $\br^r\times\{0\}$ is invariant. Then it is also invariant for the matrix $S$. Let $L_1(0)=\{l_1\,|\, (l_1,0)\in L\}$. Assume condition (a) in Proposition \ref{pri1} is satisfied, and in addition, the following hold:
\begin{enumerate}
\item
The measure $\mu_1$ is spectral with spectrum $\Lambda_1$.
\item
$S_1\Lambda_1+L_1(0)\supset\Lambda_1$, and every $\lambda_1$ in $\Lambda_1$ is an $S_1$-period for 
$$W_1(x):=\left|\frac{1}{N_1}\sum_{i=1}^{N_1}e^{2\pi i r_i\cdot x}\right|\quad(x\in\br^r).$$
\end{enumerate}
Then $\br^r\times\{0\}$ is a spectral invariant set, with spectrum 
$$\Lambda(\br^r\times\{0\}):=\bigcup_{n\geq0} (L+SL+\dots+S^{n-1}L+S^n(\Lambda_1\times\{0\})).$$
\end{corollary}

\begin{proof}
We use Theorem \ref{thi1}. We get that $\br^r\times\{0\}$ is invariant for $S$, and $l_2^0=0$ in our case. We have to check that all conditions are satisfied. Since $\br^r\times\{0\}$ is an invariant set, we have that 
$W_B(\sigma_l(x,0))=0$ for all $l\in L$ with $l_2\neq0$. This implies, using \eqref{eq1_8} that 
$$1=\sum_{l_1\in L_1(0)}W_B(\sigma_{(l_1,0)}(x,0))=\sum_{l_1\in L_1(0)}W_B(S_1^{-1}(x+l_1),0)=$$$$\sum_{l_1\in L_1(0)}\frac{1}{N^2}\sum_{i,i'=1}^{N_1}
\sum_{j=1}^{N_2(i)}\sum_{j'=1}^{N_2(i')}e^{2\pi i(r_i-r_i')\cdot(S_1^{-1}(x+l_1))}=
\sum_{l_1\in L_1(0)}\sum_{i,i'=1}^{N_1}\frac{N_2(i)}{N}\frac{N_2(i')}{N}e^{2\pi i(r_i-r_i')\cdot(S_1^{-1}(x+l_1))}$$$$
=\sum_{l_1\in L_1(0)}\left|\sum_{i=1}^{N_1}\frac{N_2(i)}{N}e^{2\pi i r_i\cdot(S_1^{-1}(x+l_1))}\right|^2=\sum_{l_1\in L_1(0)}|\hat\delta(S_1^{-1}(x+l_1))|^2,$$
where $\delta$ is the discrete measure $\delta=\sum_{i=1}^{N_1}\frac{N_2(i)}{N}\delta_{r_i}$, and $\hat\delta$ is its Fourier transform. Since $\sum_i N_2(i)=N$, the measure $\delta$ is a probability measure.

Then, the previous calculation and Proposition \ref{prs2} implies that $S_1^{-1}L_1(0)$ is a spectrum for the measure $\delta$. With \cite[Lemma 2.7]{DuJo07b} we obtain that the numbers $N_2(i)/N$ are all equal. This shows that condition (b) in Proposition \ref{pri1} is satisfied. 
Also, since $S_1^{-1}L_1(0)$ is a spectrum for $\delta$ this implies that $\# L_1(0)=N_1$ and $(A_1,\{r_1,\dots,r_{N_1}\},L_1(0))$ is a Hadamard triple. 

Finally, for $\lambda_1\in\Lambda_1$ we have $S^n(\lambda_1,0)=(S_1^n\lambda_1,0)$ so $(\lambda_1,0)$ is an $S$-period for $W_B$.

Thus all the conditions of Theorem \ref{thi1} are satisfied and the result follows.
\end{proof}

\section{Examples}\label{exam}

The following example is a spectral pair in $\br^2$ whose Hadamard triple is not reducible to $\br \times \{0\}$ , i.e., to the first coordinate. The Hadamard matrix of the system is of the form $H \otimes H$ where $H$ is the $2 \times 2$ unitary matrix of the Fourier transform on $\bz_2$, the cyclic group of order 2.

\begin{example}\label{exnr}
We give an example of an affine IFS that satisfies the Hadamard condition, but not the reducibility condition of  \cite[Theorem 3.8]{DuJo07b}.

Let $$R:=\begin{bmatrix}
4&0\\
0&4\end{bmatrix},B:=\left\{\vectr{0}{0},\vectr{0}{2},\vectr{1}{4},\vectr{1}{6}\right\}, L:=\left\{\vectr{0}{0},\vectr{2}{0},\vectr{2}{1},\vectr05\right\}$$.

Then the matrix in \eqref{eq1_5} is 
$$\frac{1}{\sqrt4}\begin{bmatrix}
1&1&1&1\\
1&1&-1&-1\\
1&-1&-1&1\\
1&-1&1&-1
\end{bmatrix}$$
which is unitary so $(R,B,L)$ is a Hadamard pair.

The subspace $\br\times\{0\}$ is invariant. Indeed, we have
$$W_B(x,y)=\left|\frac14(1+e^{2\pi i2y}+e^{2\pi i(x+4y)}+e^{2\pi i(x+6y)})\right|^2=\left|\frac14(1+e^{2\pi i2y})(1+e^{2\pi i(x+ 4y)})\right|^2.$$
And
$$\sigma_{(0,0)}(x,0)=(\frac x4,0), \sigma_{(2,0)}(x,0)=(\frac{x+2}{4},0), 
\sigma_{(2,1)}(x,0)=(\frac{x+2}{4},\frac14), \sigma_{(0,5)}(x,0)=(\frac{x}4,\frac54),$$ so 
$W_B (\sigma_{(2,1)}(x,0))=W_B(\sigma_{(0,5)}(x,0))=0$, for all $x\in\br$.

However the Hadamard triple is not reducible to $\br\times 0$ because the set $L$ does not satisfy condition \cite[Definition 3.1(iii)]{DuJo07b}: the number of vectors in $L$ that have the second component 0 is 2, and there is only one vector that has the second component 1.

\begin{theorem}\label{th4_1}
The measure $\mu=\mu_B$ is spectral.
\end{theorem}
\begin{proof}
 We use Theorem \ref{thcora}, and look for candidates for subspaces $V$ such that some finite union 
$\mathcal R=\cup_{i=1}^n(x_i,y_i)+V$ is invariant and contains some minimal invariant set $K$. Note that since $K$ is minimal, every orbit of every point in $K$ is equal to $K$. Thus we can take $(x_i,y_i)$ to be in $K$ and also be limit points of some trajectory. So we can assume $(x_i,y_i)\in X_L$, the attractor of the IFS $(\sigma_l)_{l\in L}$.

Using \eqref{eq1_3}, we have that $X_L\subset[0,\frac23]\times[0,\frac53]$.

Since $\mathcal R$ is invariant we must have for all $(x,y)\in (x_i,y_i)+V$ and $l\in L$ either $\sigma_l(x,y)\in \mathcal R$ or $W_B(\sigma_l(x,y))=0$.
Since $X_L$ is not contained in any finite union of translates of a subspace, we must have $W_B(\sigma_l(x,y))=0$ for some $(x,y)\in \mathcal R$ and some $l\in L$. 

Looking at the formula for $W_B$ we see that $W_B(x',y')=0$ iff $4y'=2k+1$ for some $k\in\bz$ or $2x'+8y'=2k+1$ for some $k\in\bz$.
Thus, the subspaces we are looking for are $V_1=\{(x,y)\,|\,y=0\}=\br\times\{0\}$ and $V_2=\{(x,y)\,|\,x+4y=0\}$.

We analyze each subspace separately, and compute what are the possible unions $\mathcal R$.

For $V_1=\br\times\{0\}$: since we have $W_B(\sigma_l(x_i,y_i))=0$ for some $l=(l_1,l_2)$ this implies 
that $4((y_i+l_2)/4)$ is of the form $2k+1$ for some $k\in\bz$. So $y_i\in\bz$. And since $(x_i,y_i)\in X_L\subset[0,2/3]\times[0,5/3]$ this implies that $y_i=0$ or $y_i=1$.

For $y_i=0$ we obtain the invariant set $\br\times\{0\}$. For $y_i=1$, since we have $\sigma_{(2,1)}(\ast,1)=(\ast,\frac12)$ and $W_B(\ast,\frac12)\neq0$, it follows that $(\ast,\frac12)\in\mathcal R$. By induction $(\ast,\frac{1}{2\cdot4^n})\in\mathcal R$, which contradicts the fact that $\mathcal R$ is a {\it finite }union of translates of $\br\times\{0\}$.

Then for $V_2=\{(x,y)\,|\, x+4y=0\}$ we will use a matrix $M$ to conjugate our IFS to another one for which $V_2$ becomes $\br\times\{0\}$, and we use Proposition \ref{prop2_8}. 
Take $M=\begin{bmatrix}
4&-1\\
0&1\end{bmatrix}$. Then $\tilde R=MRM^{-1}=R$, 
$$\tilde B=MB=B:=\left\{\vectr{0}{0},\vectr{-2}{0},\vectr{0}{1},\vectr{-2}{1}\right\},$$$$ \tilde L:=(M^T)^{-1}L=\left\{\vectr{0}{0},\vectr{0}{2},\vectr{-1}{6},\vectr{-5}{20}\right\}.$$
$$W_{\tilde B}(x,y)=\frac{1}{4}(1+e^{2\pi i (-2x)})(1+e^{2\pi iy}).$$

The subspace $V_2=\{(4t,-t)\,|\,t\in\br\}$ is mapped into $\br\times\{0\}$ by $(M^T)^{-1}$. 
We look for possible invariant unions $\mathcal R$ of the form $\cup_i(x_i,y_i)+\br\times\{0\}$.

We note that $X_{\tilde L}\subset [-\frac53,0]\times[0,\frac{20}3]$. Since $\sigma_l(\ast,y_i)=0$ for some $l\in\tilde L$ we get 
$2\cdot\frac{y_i+l_2}{4}=2k+1$ for some $k\in\bz$. Therefore $y_i$ must be an even integer in $[0,\frac{20}{3}]$. Thus $y_i\in\{0,2,4,6\}$.

We will use the following notation: if $(\ast,a)$ is in $\tilde R$ and for some $l=(l_1,l_2)$ we have that $2\cdot\frac{a+l_2}{4}$ is not an odd integer, then $W_B(\sigma_l(\ast,a))\neq0$ so $(\ast,\frac{a+l_2}4)$ must be in $\mathcal R$; we write $a\stackrel{l_2}{\rightarrow}\frac{a+l_2}4$.

If $y_i=0$, then $0\stackrel{20}{\rightarrow}\frac{20}4=5\stackrel{0}{\rightarrow}\frac54\stackrel{0}{\rightarrow}\frac5{16}\stackrel{0}{\rightarrow}\dots$, and this contradicts the finiteness of the union $\mathcal R$. 
If $y_i=2$, then $2\stackrel{2}{\rightarrow}\frac{2+2}{4}=1\stackrel{0}{\rightarrow}\frac14\stackrel{0}{\rightarrow}\dots$.
If $y_i=4$ then $4\stackrel{0}{\rightarrow}1\stackrel{0}{\rightarrow}\frac14\stackrel{0}{\rightarrow}\dots$.
If $y_i=6$ then $6\stackrel{2}{\rightarrow}2\stackrel{2}{\rightarrow}1\stackrel{0}{\rightarrow}\frac14\dots$.
In all cases we obtain a contradiction. Therefore $V_2$ does not produce invariant sets correlated to minimal invariant sets $K$ as in Theorem \ref{thcora} (one needs to apply the conjugation matrix back to $(R,B,L)$).

Thus, the only invariant sets we have to worry about are $\br\times\{0\}$ and possible $W_B$-cycles.

First, let us see what the spectrum is for $\br\times\{0\}$. We use Corollary \ref{co2}. 
We have $N_1=2$, $N_2=2$, $\{r_1,r_2\}=\{0,1\}$. The measure $\mu_1$ is associated to the IFS $\tau_0(x)=\frac x4$, $\tau_1x=\frac{x+1}4$, so it clearly has no overlap. The measure $\mu_1$ is a spectral measure with spectrum 
$$\Lambda_1:=\left\{\sum_{k=0}^n4^ka_k\,|\,a_k\in\{0,2\}\right\}.$$
(This is just a rescale of the first example of a fractal spectral measure given by Jorgensen and Pedersen in \cite{JoPe98}). 
Since $\Lambda_1\subset\bz$ we see that condition (ii) in Corollary \ref{co2} is satisfied. Hence $\br\times\{0\}$ has spectrum 
$$\Lambda(\br\times\{0\})=\bigcup_{n\geq 0}(L+SL+\dots S^{n-1}L+S^n(\Lambda_1\times\{0\})).$$

It remains to look for $W_B$-cycles. By \cite[Theorem 4.1]{DuJo07a} we have to compute the lattice
$\Gamma=\{\gamma\in\br^2\,|\,\gamma\cdot b\in\bz\mbox{ for all }b\in B\}$. So if $(x,y)\in \Gamma$ then $2x, x+4y,x+6y\in\bz$ so
$2y\in\bz$ and $x\in\bz$, which means that $\Gamma=\bz\times\frac12\bz$. To find the $W_B$-cycles we intersect $\Gamma$ with $X_L$. Since $X_L\subset[0,\frac23]\times[0,\frac53]$ it follows that the only candidates for $W_B$-cycle points are 
$(0,0),(0,\frac12),(0,1)$. 

$(0,0)$ is the trivial cycle, but we can discard it because it is contained in $\br\times\{0\}$. It is easy to check that the other two are not $W_B$-cycles. 

We conclude that the invariant set $\br\times\{0\}$ contains all minimal invariant sets, and therefore, by Proposition \ref{propsupp}, the spectrum 
of $\mu$ is $\Lambda(\br\times\{0\})$.
\end{proof}
\end{example}
\begin{acknowledgements}
The authors are pleased to thank the following for helpful discussions at various times:  John Benedetto, Deguang Han, Keri Kornelson, Paul Muhly, Erin Pearse, Steen Pedersen, Gabriel Picioroaga, Karen Shuman, Bob Strichartz, Qyiu Sun and Yang Wang.
\end{acknowledgements}

\bibliographystyle{alpha}	
\bibliography{spec2}

\end{document}